\documentclass[11pt]{article}
\usepackage{amssymb,amsmath,bm,graphicx}
\usepackage{mathrsfs}
\usepackage{constants}
\usepackage{hyperref}
\usepackage{enumitem}
\usepackage{color}

\date{}
\begin{document}
\newcommand{\bea}{\begin{eqnarray}}
\newcommand{\ena}{\end{eqnarray}}
\newcommand{\beas}{\begin{eqnarray*}}
\newcommand{\enas}{\end{eqnarray*}}
\newcommand{\beq}{\begin{equation}}
\newcommand{\enq}{\end{equation}}
\def\qed{\hfill \mbox{\rule{0.5em}{0.5em}}}
\newcommand{\bbox}{\hfill $\Box$}
\newcommand{\ignore}[1]{}
\newcommand{\ignorex}[1]{#1}
\newcommand{\wtilde}[1]{\widetilde{#1}}
\newcommand{\qmq}[1]{\quad\mbox{#1}\quad}
\newcommand{\qm}[1]{\quad\mbox{#1}}
\newcommand{\nn}{\nonumber}
\newcommand{\Bvert}{\left\vert\vphantom{\frac{1}{1}}\right.}
\newcommand{\To}{\rightarrow}
\newcommand{\E}{\mathbb{E}}
\newcommand{\Var}{\mathrm{Var}}
\newcommand{\Cov}{\mathrm{Cov}}
\newcommand{\Corr}{\mathrm{Corr}}
\newcommand{\dist}{\mathrm{dist}}
\newcommand{\diam}{\mathrm{diam}}
\makeatletter
\newsavebox\myboxA
\newsavebox\myboxB
\newlength\mylenA
\newcommand*\xoverline[2][0.70]{%
    \sbox{\myboxA}{$\m@th#2$}%
    \setbox\myboxB\null
    \ht\myboxB=\ht\myboxA%
    \dp\myboxB=\dp\myboxA%
    \wd\myboxB=#1\wd\myboxA
    \sbox\myboxB{$\m@th\overline{\copy\myboxB}$}
    \setlength\mylenA{\the\wd\myboxA}
    \addtolength\mylenA{-\the\wd\myboxB}%
    \ifdim\wd\myboxB<\wd\myboxA%
       \rlap{\hskip 0.5\mylenA\usebox\myboxB}{\usebox\myboxA}%
    \else
        \hskip -0.5\mylenA\rlap{\usebox\myboxA}{\hskip 0.5\mylenA\usebox\myboxB}%
    \fi}
\makeatother

\newtheorem{theorem}{Theorem}[section]
\newtheorem{corollary}[theorem]{Corollary}
\newtheorem{conjecture}[theorem]{Conjecture}
\newtheorem{proposition}[theorem]{Proposition}
\newtheorem{lemma}[theorem]{Lemma}
\newtheorem{definition}[theorem]{Definition}
\newtheorem{example}[theorem]{Example}
\newtheorem{remark}[theorem]{Remark}
\newtheorem{case}{Case}[section]
\newtheorem{condition}{Condition}[section]
\newcommand{\proof}{\noindent {\it Proof:} }

\title{{\bf\Large Normal approximation for associated point processes via Stein's method with applications to determinantal point processes}}
\author{Nathakhun Wiroonsri \thanks{Department of Mathematics, King Mongkut's University of Technology Thonburi, Bangkok 10140, Thailand. Email: nathakhun.wir@mail.kmutt.ac.th }  \\ King Mongkut's University of Technology Thonburi}

\maketitle

\begin{abstract} 
We use Stein's method to provide non asymptotic $L^1$ bounds (also known as Wasserstein bounds) to the normal for functionals of associated point processes. As for supporting tools, we use the connection between association and $\alpha$-mixing properties that was recently uncovered by \cite{PDL17}. We apply our main results to determinantal point processes which are known to be negatively associated. A potential application to point processes in the Laguerre-Gaussian family is also presented.
\end{abstract}

\section{Introduction}

In this work, we obtain $L^1$ bounds for functionals of determinantal point processes which are generalized from the CLT results without rates of convergence in \cite{PDL17}. We recall that the $L^1$ distance between the distributions $\mathcal{L}(X)$ and $\mathcal{L}(Y)$ of real valued random variables $X$ and $Y$ is given by
\bea \label{Wasdef}
d_1 \big(\mathcal{L}(X),\mathcal{L}(Y)\big) &=& \int_{-\infty}^\infty |P(X \le t) - P(Y \le t)| dt \nn \\
                                      &=& \sup_{h \in \mathcal{H}_1} |\E h(X)-\E h(Y)| ,
\ena 
where $\mathcal{H}_1 =\{h: |h(y)-h(x)| \le |y-x|\}$. The $L^1$ distance is also known as the Wasserstein distance.
Some of the readers might be more familiar with the $L^\infty$ or Kolmogorov distance defined by $\sup_{t \in \mathbb{R}}|P(X \le t) - P(Y \le t)|$. In general the two distances are not comparable as it is known that $\sup_{x \in \mathbb{R}} |P(X \le x) - P(Y \le x)| \le \sqrt{2cd_1\big({\cal L}(X),{\cal L}(Y)\big)}$ for some $c>0$ but not conversely (See Proposition 1.2 of \cite{Ross11}).

Point processes (sometimes called random point fields) are one of the most powerful tools in probability that have been widely developed from their root, Poisson point process, and have been applied to sciences and technology ranging from applied mathematics to astronomy. 
Determinantal point processes, detailed in Section \ref{sec:app}, are well known examples of point processes that are characterized by determinants of kernel functions and known to be models of repulsion. Although they first appeared in quantum physics (see \cite{VY01} for example), nowadays they have been developed in many branches of applied sciences and mathematics such as random matrix theory (\cite{HKPV09}), statistics (\cite{BH16}), wireless communication network (\cite{MS14}, \cite{DZH15}) and machine learning (\cite{KT12}).

In this work, we have two main interests. First, we aim to obtain $L^1$ bounds to the normal for the sum of some functions of an associated point process. Secondly, we apply the bounds to functionals of a determinantal point process. We postpone the full detail of the models to the remaining sections and discuss here some known results related to our work. Though we are unaware of any results that provide a rate of convergence under the setting of this work, there are several papers that obtain ones for Poisson point processes. As it is not our main interest, let us refer to one of the most recent and most related paper to our work (\cite{RS13}). For the reader who is interested in more details towards that direction, please see the references therein.  
For $\eta$ a Poisson point process on a Borel space $X$, the work \cite{RS13} obtained an $L^1$ bound to the normal of a $U$-statistic of $\eta$ defined by 
\beas
F(\eta) = \sum_{\left\{x_1,\ldots x_k\right\} \in \eta^k_{\neq}}f(x_1,\ldots,x_k)
\enas 
where $\eta^k_{\neq}$ is the set of all $k$-tuples of distinct points of $\eta$ and $f$ is an integrable function on $X^k$. The bound is given as
\beas
d_1\left(\mathcal{L}\left(\frac{F-\E F}{\sqrt{\Var(F)}}\right),\mathcal{L}(Z)\right) \le 2k^{7/2} \sum_{1\le i \le j \le k} \frac{\sqrt{M_{i,j}(F)}}{\Var(F)}
\enas
where $M_{i,j}(F)$ are sums of certain fourth moment integrals and $Z$ is a standard normal random variable. The result can be specialized to the case when $f$ is independent of the intensity parameter $\lambda$ of the Poisson process where the authors obtained the $L^1$ bound 
\beas
d_1\left(\mathcal{L}\left(\frac{F-\E F}{\sqrt{\Var(F)}}\right),\mathcal{L}(Z)\right)   \le  C_f \lambda^{-1/2},
\enas 
with a constant $C_f$ depending on $f$. 

For comparison, in this work we obtain an $L^1$ bound to the normal of order $n^{-d/(4d+2)}$ between a functional of a determinantal point process on $\mathbb{R}^d$ with $d \in \mathbb{N}_1$ defined later in \eqref{fdef} and the normal under a few assumptions specified in Section \ref{sec:app}. Our bound does not depend on an intensity parameter but $n^d$ is the volume of the domain that we consider. To be more precise, the work \cite{RS13} considered the process on the entire space $X$ and obtained a bound based on the intensity parameter $\lambda$, on the other hand, this work obtains a bound based on the size $n^d$ of a subset of the entire space $\mathbb{R}^d$. The functional has a similar form as the $U$-statistic above but we require $f$ to be bounded and do not fix $k$. Some known results related to this functional will also be mentioned after the definitions are stated properly. As for example we consider the Laguerre-Gaussian family of determinantal point processes (\cite{BL16}).

Determinantal point processes are well known to be negatively associated. As a result, we proceed to our goal by first providing $L^1$ bounds to the normal for the sum of functions of associated point processes in general. Then we write the functional of a determinantal point process as a sum of the same form. Positive association has been found frequently in statistical physics models and has been used to develop normal approximations for random fields under different methods, see for instance \cite{New80}, \cite{New83}, \cite{CG84}, \cite{Bir88}, \cite{Bul95} and \cite{GW18}. Negative association has also been used in the same manner and has well known applications related to permutation distributions, see for example \cite{JP83}, \cite{CW09} and \cite{Wir18}. However, association properties had not appeared in the context of normal approximation for point processes until \cite{PDL17} obtained a CLT for their functionals under some certain conditions. As we provide rates of convergence in this work, we assume stronger conditions than \cite{PDL17}.

The main tool that we use in this work is Stein's method, introduced by \cite{Stein72}, which is now one of the most powerful methods to prove convergence in distribution as it has main advantages that it provides non-asymptotic bounds on the distance between distributions, and that it can handle various situations involving dependence. Thus far, many applications in several areas such as statistics, statistical physics and applied sciences have been developed using this method. For more details about the method in general, see the text \cite{CGS11} and the introductory notes \cite{Ross11}.  Though we are unaware of the use of Stein's method on normal approximation for a point process, it was used earlier to obtain bounds between point processes under various assumptions and Poisson processes. The following is an incomplete list of references on Poisson process approximation using Stein's method, \cite{Bar88}, \cite{BB92}, \cite{CX04}, \cite{CX06} and \cite{CX11}. The method was also used to bound distances between sums of positively and negatively associated random variables and some well-known targeted distributions (see \cite{GW18} and \cite{Wir18} for the normal and \cite{Dal13} for the Poisson).

To support Stein's method, we definitely use the association property and an additional condition called $\alpha$-mixing (also known as strong mixing) for point processes whose connection was recently uncovered by \cite{PDL17}. The $\alpha$-mixing condition, introduced by \cite{Ros56}, is a useful tool in probability theory used to measure dependency between two $\sigma$-algebras. It has many probabilistic properties and applications which are even covered in the whole book \cite{Dou94} (see also the survey \cite{Bra05}). The $\alpha$-mixing conditions alone have been used broadly to prove CLTs for dependent random variables (see \cite{PDL17} and the references therein).

The remaining of this work is organized as follows. We first recall some necessary definitions, notations and previous results in Section \ref{sec:def}. In Section \ref{sec:main} we provide and prove $L^1$ bounds for a sum of functions of associated point processes in general. Then we devote Section \ref{sec:app} to discuss an application to determinantal point processes.

\section{Some background and definitions} \label{sec:def}

In this section we recall some definitions and state some notations. We first recall the definition of point processes on $\mathbb{R}^d$ in general. In this paper, unless otherwise stated, $X$ represents a point process and $x$ denotes its realization. Letting $\text{card}(A)$ be the cardinality of $A \subset \mathbb{R}^d$, $x\subset \mathbb{R}^d$ is said to be \textit{locally finite} if $\text{card}(x \cap B) < \infty$ for all bounded sets $B \subset \mathbb{R}^d$ and for
\beas
\Omega = \left\{x \subset \mathbb{R}^d:\text{card}(x \cap B)<\infty,\forall \text{  bounded  } B\subset \mathbb{R}^d\right\},
\enas
elements of $\Omega$ are called \textit{locally finite point configurations}. For a point process $X \subset \mathbb{R}^d$, we denote 
\beas
N(A) = \text{card}(X \cap A)
\enas 
and recall that $X$ is said to be \textit{simple} if for any ${\bf a} \in \mathbb{R}^d$, $N(\{{\bf a}\}) \in \{0,1\}$ a.s. and \textit{locally finite} if it takes values in $\Omega$. In this work, we focus only on locally finite simple point processes on $\mathbb{R}^d$. Here, with $\mathbb{Z}$ the set of integers and $k \in \mathbb{Z}$, we let $\mathbb{N}_k=[k,\infty) \cap \mathbb{Z}$. With $\mathscr{B}(A)$ the Borel $\sigma$-algebra of $A \subset \mathbb{R}^d$, we equip $\Omega$ with,
\beas
\mathcal{E}(A) = \sigma\left(\{x \in \Omega,\text{card}(x \cap B)=m\}:B \in \mathscr{B}(A),m \in \mathbb{N}_0\right),
\enas 
the smallest $\sigma$-algebra generated by $\{x \in \Omega,\text{card}(x \cap B)=m\},B \in \mathscr{B}(A),m \in \mathbb{N}_0$. A point process can be determined by one of the three characterizations: its finite dimensional distributions, its void probabilities and its generating functionals. However, we do not focus on this theoretical background. If interested, the reader may check out the books \cite{DV03}, \cite{MW04} and \cite{IPSS08}.

Next we recall that association properties for point processes are defined as follows.
A point process $X$ is said to be \textit{negatively associated} if for all families of pairwise disjoint Borel sets $(A_i)_{1\le i \le k}$ and $(B_j)_{1\le j \le l}$ such that
\bea \label{disj}
(\cup_iA_i) \cap (\cup_jB_j) = \emptyset,
\ena
and for all coordinate wise increasing functions $\psi:\mathbb{N}^k\rightarrow \mathbb{R}$ and $\phi:\mathbb{N}^l\rightarrow \mathbb{R}$,
\begin{multline*}
\E\left[\psi\left(N(A_1),\ldots,N(A_k)\right)\phi\left(N(B_1),\ldots,N(B_l)\right)\right] \\
 \le \E\left[\psi\left(N(A_1),\ldots,N(A_k)\right)\right] \E \left[\phi\left(N(B_1),\ldots,N(B_l)\right)\right] .
\end{multline*}
Similarly, a point process is said to be \textit{positively associated} if it satisfies the reverse inequality for all families of pairwise disjoint Borel sets $(A_i)_{1\le i \le k}$ and $(B_j)_{1\le j \le l}$ but not necessarily satisfying \eqref{disj}. In addition, a point process is said to be \textit{associated} if it is either negatively or positively associated. In some papers, association only refers for short to positive association (see \cite{Bir88} for example).

Before moving on to the next section, we provide a short review of Stein's method. We first recall that, for real valued random variables $X$ and $Y$ and some family of functions $\mathcal{H}$, the distance between the distributions of $X$ and $Y$ can be defined of the form
\bea \label{disdef}
d_\mathcal{H} (\mathcal{L}(X),\mathcal{L}(Y)) = \sup_{h \in \mathcal{H}} |\E h(X)-\E h(Y)| .
\ena 
In this work we consider the $L^1$ distance which corresponds to the case where $\mathcal{H}=\mathcal{H}_1$ as defined in \eqref{Wasdef}.
Stein's method was motivated from the fact that $W$ has the standard normal distribution, denoted ${\cal N}(0,1)$, if and only if 
\beas
\E W f(W) = \E f'(W)
\enas
for all absolutely continuous functions $f$ with $\E |f'(W)| < \infty$. This equation and the distance in \eqref{disdef} lead to the differential equation
\bea \label{steineq}
h(w) - \E h(Z) = f'_h(w)-wf_h(w)
\ena
where $Z \sim \mathcal{N}(0,1)$ and $h \in \mathcal{H}$. Taking the supremum over all $h \in \mathcal{H}$ to the expectation on the left hand side of \eqref{steineq} with $w$ replaced by a variable $W$ yields the distance between $W$ and $Z$ in \eqref{disdef}. Thus one can instead handle the expectation on the right hand side using the bounded solution $f_h$ of \eqref{steineq} for the given $h$. Using this device, Stein's method has uncovered an alternative way to show convergence in distribution with additional information on the finite sample distance between distributions. 

A property of $f_h$, the solution of \eqref{steineq} when $h \in \mathcal{H}_1$, is stated in the lemma below. For a real valued function $\varphi(u)$ defined on the domain $\mathcal{D}$, let $|\varphi|_{\infty}= \sup_{x \in \mathcal{D}} |\varphi(x)|$. We include in this definition the $|\cdot|_\infty$ norm of vectors and matrices, for instance, by considering them as real valued functions of their indices. 
\begin{lemma} [\cite{CGS11}] \label{Waslem}
Let $h \in \mathcal{H}_1$ where $\mathcal{H}_1$ is given in \eqref{Wasdef} and $f_h$ be the bounded solution of \eqref{steineq}. Then
\bea \label{steineq:bounds}
|f_h|_{\infty} \le 2, |f_h'|_{\infty} \le \sqrt{\frac{2}{\pi}} \qmq{and} |f_h''|_{\infty} \le 2.
\ena
\end{lemma}

\section{$L^1$ bounds for associated point processes} \label{sec:main}

In this section, we generalize Section 3 of \cite{PDL17} that proved a CLT under the same setting. Let $X \in \Omega$ be an associated point process and
\bea \label{Ydef}
Y_{\bf i} = f_{\bf i}(X \cap C_{\bf i}) - \E f_{\bf i}(X \cap C_{\bf i}), {\bf i} \in \mathbb{Z}^d
\ena
where $f_{\bf i}:\Omega \rightarrow \mathbb{R}$ are real valued measurable functions, $C_{\bf i}$, ${\bf i} \in \mathbb{Z}^d$ are defined as the $d$-dimensional unit cube centered at ${\bf i}$. Note that the union of $C_{\bf i}$ forms a covering of $\mathbb{R}^d$. 

We let ${\bf 1} \in \mathbb{Z}^d$ denotes the vector with all components $1$, and write inequalities such as ${\bf a} < {\bf b}$ for vectors ${\bf a} , {\bf b} \in \mathbb{R}^d$  when they hold componentwise. In this work, we consider 
\bea \label{blockdef}
S_{{\bf k}}^n =  \sum_{{\bf i} \in B_{{\bf k}}^n} Y_{{\bf i}} \qmq{where}
B_{{\bf k}}^n = \left\{ {\bf i} \in \mathbb{Z}^d: {\bf k} \le {\bf i} < {\bf k}+n{\bf 1} \right\}.
\ena   
The work \cite{PDL17} obtained a CLT for the sum above with $B_{{\bf k}}^n$ replaced by any sequences of strictly increasing finite domains of $\mathbb{Z}^d$ and $C_{\bf i}$ be any $d$-dimensional cube centered at $x_{\bf i} = R \cdot {\bf i}$ with fixed $R>0$ and with fixed side length $s \ge R$. 

Although we restrict our interest to the case that the sum is over a square block and that $C_{\bf i}$ is the $d$-dimensional unit cube centered at ${\bf i}$, we will remark right after our main theorem that generalizing to a cube with fixed $R$ and $s$ only affects our results by a constant factor and therefore it suffices to study only the case that $R=s=1$. Also, we only consider the square block sum of size $n^d$ but we remark that our method works for some more specific types of sum as well but with more complicated computation. 

Prior to stating the main theorem, we state the definition of the nth order intensity functions of point processes with respect to Lebesque measure. 
Let $n \in \mathbb{N}_1$ and $X \in \Omega$. If there exists a non negative function $\rho_n : \left(\mathbb{R}^d\right)^n \rightarrow \mathbb{R}$ such that
\beas
\E\left[\sum_{\substack{{\bf x}_1,\ldots,{\bf x}_n \in X \\ \text{all distinct}}}f({\bf x}_1,\ldots,{\bf x}_n)\right] 
= \int_{\left(\mathbb{R}^d\right)^n} f({\bf x}_1,\ldots,{\bf x}_n) \rho_n({\bf x}_1,\ldots,{\bf x}_n) d{\bf x}_1 \ldots d{\bf x}_n,
\enas
for all locally integrable functions $f : \left(\mathbb{R}^d\right)^n \rightarrow \mathbb{R}$, then $\rho_n$ is called the \textit{nth order intensity function} of $X$.
Now for ${\bf x},{\bf y} \in \mathbb{R}^d$, let
\bea \label{ddef}
D({\bf x},{\bf y}) = \rho_2({\bf x},{\bf y})-\rho_1({\bf x})\rho_1({\bf y}).
\ena
It follows that
\beas
\Cov(N(A),N(B)) = \int_{A \times B} D({\bf x},{\bf y}) d{\bf x} d{\bf y}.
\enas

The following theorem provides an $L^1$ bound of order $n^{-d/(4d+2)}$ between the standardized $S_{{\bf k}}^n$ and the normal. In the following, we use the notation $\left\|Y\right\|_p = (\E|Y|^p)^{1/p}$ for $p \in \mathbb{N}_1$.

\begin{theorem} \label{mainthm1}
For $d \in \mathbb{N}_1$, let $X$ be an associated point process on $\mathbb{R}^d$, $S_{{\bf k}}^n$ be as in \eqref{blockdef} with $Y_{\bf i}$ given in \eqref{Ydef} with $R=s=1$, ${\bf k} \in \mathbb{Z}^d$ and $\sigma_{n,{\bf k}}^2 = \Var(S_{{\bf k}}^n)$. Assume that the following conditions are satisfied:
\begin{enumerate}[label=(\alph*)]
	\item The first two intensity functions of $X$ are well defined; \label{a}
	\item $\sup_{{\bf i} \in \mathbb{Z}^d} \left\|Y_{\bf i}\right\|_3 = M < \infty$; \label{b}
	\item $\sup_{|{\bf x}-{\bf y}|_\infty \ge r} D({\bf x},{\bf y}) \le \kappa e^{-\lambda r}$ for some $\kappa,\lambda > 0$; \label{c}
	\item $\sigma_{n,{\bf k}}^2 \ge \gamma n^d$ for some $\gamma>0$. \label{d}
\end{enumerate}
Then, with $Z$ be a standard normal random variable,
\bea \label{mainbound}
d_1\left(\mathcal{L}\left(\frac{S_{{\bf k}}^n}{\sigma_{n,{\bf k}}}\right),\mathcal{L}(Z)\right) &\le&
\frac{C_{1,d,M,\kappa,\gamma}}{n^{d/(4d+2)}} + \frac{C_{2,d,M,\kappa,\gamma} n^{d(4d+1)/(6d+3)}}{\exp\left(\theta_{d,M,\kappa,\gamma}n^{d/(4d+2)}\right)} \nn\\
       && \hspace{30pt}+ \frac{C_{3,d,M,\kappa,\gamma} n^{7d/6}}{\exp\left(2\theta_{d,M,\kappa,\gamma}n^{d/(4d+2)}\right)},
\ena
where 
\bea \label{thetadef}
\theta_{d,M,\kappa,\gamma} =\frac{\lambda}{3}\left( \frac{\sqrt{2\gamma}\kappa^{1/3}\left((4\mu_{\lambda}+2\nu_{\lambda})^d-\left(2\nu_{\lambda}\right)^d\right)}{18^{d+1}\sqrt{\pi}dM}\right)^{1/(2d+1)},
\ena
\beas
&&C_{1,d,M,\kappa,\gamma} = \left(\frac{9\cdot 36^d M^{4d+3}\left((4\mu_{\lambda}+2\nu_{\lambda})^d-\left(2\nu_{\lambda}\right)^d\right)^{2d} }{\gamma^{2d+(3/2)}\pi^d}\right)^{1/(2d+1)} \\ 
   && \hspace{100pt} \times \left(\frac{1}{(2d)^{2d/(2d+1)}}+2(2d)^{1/(2d+1)}\right) ,
\enas
\bea \label{Cdef}
C_{2,d,M,\kappa,\gamma} = \frac{3\cdot 6^d \kappa^{1/3} M^2 \theta_{d,M,\kappa,\gamma}^{4d/3}}{\sqrt{\pi}\gamma}, \ \ \ 
C_{3,d,M,\kappa,\gamma} = \frac{2^{d+1} \kappa^{2/3} M}{\sqrt{\gamma}},
\ena
and
\bea \label{munudef}
\mu_\lambda = \frac{e^{\frac{2\lambda}{3}}}{\left(e^{\frac{\lambda}{3}}-1\right)^2}, \ \ \ \nu_\lambda = \frac{e^{\lambda}}{\left(e^{\frac{\lambda}{3}}-1\right)^2}.
\ena
\end{theorem}

\begin{remark}
The only term in our main bound \eqref{mainbound} that contributes the rate of $n^{-d/(4d+2)}$ is the first term since the last two terms decay exponentially in $n$. 
\end{remark}

\begin{remark}
We remark that our assumption \ref{c} is stronger than the one in \cite{PDL17} that only requires a decaying rate $o(r^{-(3d+\epsilon)})$ for some $\epsilon>0$. However, rates of convergence for the CLT were not provided there. 
\end{remark}

\begin{remark} \label{nsnd}
By following the same proof, our assumption \ref{d} can be relaxed by replacing $n^d$ by $n^s$ with $s > \left(\frac{4d+2}{4d+3}\right)d$. However, the rate of $n^{-d/(4d+2)}$ will be replaced by a slower rate of $n^{-\left(\frac{4d+3}{4d+2}\right)s+d}$. We will see in the proof below that the only terms that contribute to the rate of convergence are \eqref{1st} and \eqref{2nd} and by replacing $n^d$ by $n^s$ in $\sigma_{n,{\bf k}}^2 \ge \gamma n^d$, these two terms will have orders $n^{d-s}/l$ and $l^{2d}/n^{3s/2-d}$, respectively. Setting $l = O(n^r)$ for some $r>0$, the optimal value of $r$ is $s/(4d+2)$. Plugging this back into  $n^{d-s}/l$, we have that this term is of order $n^{-\left(\frac{4d+3}{4d+2}\right)s+d}$. Since we need $\left(\frac{4d+3}{4d+2}\right)s-d>0$, it follows that $s > \left(\frac{4d+2}{4d+3}\right)d$.
\end{remark}

\begin{remark}
The same result as in Theorem \ref{mainthm1} still holds if one replaces $B_{\bf k}^n$ ,$n \in \mathbb{N}_1$ by any increasing sequences of indexes $I_n$ of the same sizes as $B_{\bf k}^n$, for example one may let $I_n = \{{\bf i} \in \mathbb{Z}^d:{\bf i} \in [-n/2,n/2]^d\}$. Also, replacing $B_{\bf k}^n$ by any increasing hypercubes of any sizes will only affect the constant. Our proving technique also works for some more specific types of $I_n$ with more complex computations, e.g. rectangular or any types of increasing sequence of indexes that cover hypercubes of side length $\alpha n$ and are covered by hypercubes of side length $\beta n$ for some fixed $0<\alpha<\beta$.
\end{remark}

\begin{remark}
Replacing the unit cubes $C_{\bf i}$ by $d$-dimensional cubes centered at $x_{\bf i} = R \cdot {\bf i}$ with fixed $R>0$ and with fixed side length $s$ where $s=R > 1$ in \eqref{Ydef} obviously does not affect the proof below and the result only changes by a constant factor. Letting $s>R>0$ affects the terms \eqref{2nd}, \eqref{3rdi1}, \eqref{4thj1} and \eqref{5th} in the proof of Theorem \ref{mainthm1} below. However, as $s$ and $R$ are fixed, it does not change the rate of convergence in the main result. For example, if $s>R=1$, the second line of \eqref{3rdi1} will be replaced by
\beas
\frac{1}{\sigma^2} \sqrt{\frac{2}{\pi}}l^d(l+2(s-1))^dM^2 \sum_{\bf i}\alpha^{1/3}_{(l+2(s-1))^d,n^d-(3l-2(s-1))^d}(l-2(s-1)), 
\enas   
which only changes the numerical value of the constant.
\end{remark}

To prove the theorem, we first introduce some notations. Note that if we assume the conditions \ref{a}-\ref{d} in the statement of Theorem \ref{mainthm1}, the proof will be the same for any fixed ${\bf k} \in \mathbb{Z}^d$. Hence we only consider $S_{\bf 1}^n$ and for simplicity in the proof we denote $S_n = S_{\bf 1}^n$, $B_n = B_{\bf 1}^n$, $\sigma_n^2=\sigma_{n,{\bf 1}}^2$ and let 
\bea \label{Wdef}
W_n = S_n/\sigma_n.
\ena

We use the technique in Section 2.1 of \cite{GW18} decomposing the sum $S_n$ over the block $B_n$
into sums over smaller, disjoint blocks whose side lengths are at most some integer $l$. That is, 
for $1 \le l \le n$, we uniquely write $n=(m-1)l+r$ with $m \ge 1$ and $1 \le r \le l$ and correspondingly decompose $B_n$ into $m^d$ disjoint blocks $D_{\bf i}^l, {\bf i} \in [m]^d$, where there are $(m-1)^d$ `main' blocks having all sides of length $l$, and $m^d-(m-1)^d$ remainder blocks having all sides of length $r$ or $l$, with at least one side of length $r$.

To be more precise, 
\beas
D_{\bf i}^l&=&\big\{{\bf j} \in \mathbb{Z}^d: (i_s-1)l+1 \le j_s \le i_sl \text{  for  } i_s \not =m, \\
                                                   &&\hspace{50pt}(m-1)l+1 \le j_s \le (m-1)l+r \text{  for  } i_s =m \big\}.
\enas
It is easy to see that for ${\bf i} \in [m-1]^d$, the vectors indexing the `main blocks', we have 
\bea \label{Di.block}
D_{\bf i}^l = B_{({\bf i}-{\bf 1})l+{\bf 1}}^l  \text{ \ \ for \ \ } {\bf i} \in [m-1]^d,
\ena
and if $r=l$ then $D_{{\bf i}}^l$ is given by \eqref{Di.block} for all ${\bf i} \in [m]^d$. Furthermore, it is straightforward to verify that the elements of the collection $\{D_{\bf i}^l, {\bf i} \in [m]^d\}$ is a partition of $B_n$.

Letting
\bea \label{xidef}
\xi_{\bf i}^{l} = \sum_{{\bf t} \in D_{\bf i}^l}Y_{\bf t} \text{ \ \ for \ \  } {\bf i} \in [m]^d,
\ena 
we see that $W_n$ can also be writen as
\beas
W_n = \sum_{{\bf i} \in [m]^d} \frac{\xi_{\bf i}^{l}}{\sigma_n} 
\enas
which has mean zero and variance one. 

In the proof, we use the following two notations
\bea \label{Whdef}
W_{{\bf i},n} = \sum_{\substack{{\bf k} \in [m]^d \\ |{\bf k} - {\bf i}|_\infty \le 1}}\frac{\xi_{\bf k}^l}{\sigma_n} 
\text{, \ \ \ } W_{{\bf i},n}^*=W_n-W_{{\bf i},n}, \nn \\
W_{{\bf i},{\bf j},n} = \sum_{\substack{{\bf k} \in [m]^d \\ |{\bf k} - {\bf i}|_\infty \wedge |{\bf k} - {\bf j}|_\infty \le 1}}\frac{\xi_{\bf k}^l}{\sigma_n} 
\text{ \ \ and \ \  } W_{{\bf i},{\bf j},n}^*=W_n-W_{{\bf i},{\bf j},n}
\ena
for ${\bf i},{\bf j} \in [m]^d$ with $a \wedge b = \min(a,b)$. Note that the ones in the upper line denote the sum of $\xi_{\bf k}^l$ over all ${\bf k}$ that are within distance $1$ from ${\bf i}$ and its complement, respectively, and the ones in the lower line denote the similar variables with ${\bf i}, {\bf j}$ replacing ${\bf i}$. 

Next we state the definition of the $\alpha$-mixing coefficient. For a probability space $(\mathcal{X},\mathcal{F},P)$ and $\mathcal{A},\mathcal{B}$ be two sub $\sigma$-algebras of $\mathcal{F}$, the \textit{$\alpha$-mixing coefficient} between $\mathcal{A}$ and $\mathcal{B}$ is defined as,
\beas
\alpha(\mathcal{A},\mathcal{B}) = \sup\left\{|P(A\cap B) - P(A)P(B)|:A \in \mathcal{A}, B \in \mathcal{B}\right\}.
\enas
It is clear that the $\alpha$-mixing coefficient can be viewed as a measurement of dependence between two sub $\sigma$-algebras under the same probability space. Particularly, it is zero if and only if $\mathcal{A}$ and $\mathcal{B}$ are independent. For a point process $X \subset \mathbb{R}^d$, $a,b,c \ge 0$, let
\bea \label{alphadef2}
\alpha_{a,b}(c) = \sup\left\{\alpha(\mathcal{E}(A),\mathcal{E}(B)):|A| \le a, |B| \le b, \dist(A,B) \ge c\right\}
\ena
with the convention that $\alpha_{a,\infty}(c) = \sup_b \alpha_{a,b}(c)$ where $\dist(A,B) = \inf_{{\bf x} \in A,{\bf y}\in B}|{\bf y}-{\bf x}|_1$ and $|{\bf x}|_1=\sum_{i=1}^d |x_i|$ is the $L^1$ vector norm of ${\bf x} = (x_1,\ldots,x_d)\in \mathbb{R}^d$.

Now we state the following two lemmas used in the proof that follows. The first one from \cite{Rio93} bounds the covariance between two random variables by their norms and the $\alpha$-mixing coefficient.
\begin{lemma} [\cite{Rio93}] \label{lemrio}
Let $X$ and $Y$ be random variables on the same probability space and measurable with respect to $\mathcal{A}$ and $\mathcal{B}$, respectively. Then
\beas 
|\Cov(X,Y)| \le \alpha^{1/r}(\mathcal{A},\mathcal{B})\left\|X\right\|_p \left\|Y\right\|_q, 
\enas
for all $p,q,r \in [1,\infty]$ such that $1/p+1/q+1/r = 1$.
\end{lemma} 
Next lemma proved in \cite{PDL17} provides bounds for $\alpha_{a,b}(c)$, defined in \eqref{alphadef2}, in term of $D(\cdot,\cdot)$ as given in \eqref{ddef}.
\begin{lemma} [\cite{PDL17}] \label{lempdl}
Let $X$ be an associated point process on $\mathbb{R}^d$ whose first two intensity functions are well defined. Then for all $a,b>0$ and $c \ge 0$,
\beas
\alpha_{a,b}(c) \le ab \sup_{|{\bf x}-{\bf y}|\ge c} |D({\bf x},{\bf y})| \text{ \ \  \ \ \ \ \ and }
\enas
\beas
\alpha_{a,\infty}(c) \le as_d \int_c^\infty t^{d-1} \sup_{|{\bf x}-{\bf y}| = t} |D({\bf x},{\bf y})| dt.
\enas
\end{lemma}

In the following we will apply the identities
\bea \label{Sum.krk4}
\sum_{k=1}^{n-1} (n-k)w^{k}= \frac{w\left((n-1)-nw+w^n\right)}{(w-1)^2} \text{ \ \ for \ \ } w \ne 1,
\ena
and
\bea \label{Sum.krk3}
n + \sum_{a=1}^{n-1} (n-a)(v^a+v^{-a})= \frac{v^{1-n} \left(v^n-1\right)^2}{(v-1)^2} \text{ \ \ for \ \ } v \ne 1.
\ena
Also for simplicity in the proof we will drop $l$ in $\xi_{\bf i}^l$ and $n$ in $W_n$, $W_{{\bf i},n}$ and $W_{{\bf i},{\bf j},n}$ given in \eqref{xidef} and \eqref{Whdef}, respectively.

\bigskip

\noindent {\textbf{Proof of Theorem \ref{mainthm1}:}}  
Since we prove non-asymptotic bounds for any fixed $n \in \mathbb{N}_1$, we drop $n$ in all the notations in the proof for simplicity. For given $h \in \mathcal{H}_1$ let $f$ be the unique bounded solution to \eqref{steineq}. Then, by Lemma \ref{Waslem}, we have
\bea \label{steineq:bounds}
|f'|_{\infty} \le \sqrt{\frac{2}{\pi}} \qmq{and} |f''|_{\infty} \le 2.
\ena
Below, unless otherwise stated, all indexes under $\sum$ are taken over $[m]^d$. Denoting $\sigma_{\bf i}^2 = \Var(\xi_{\bf i})$ and $\sigma_{{\bf i},{\bf j}}=\Cov(\xi_{\bf i},\xi_{\bf j})$ for ${\bf i} \ne {\bf j} \in [m]^d$, using that $\sigma^2 = \sum_{\bf i} \sigma^2_{\bf i} + \sum_{{\bf i} \ne {\bf j}} \sigma_{{\bf i},{\bf j}}$, with $W_{\bf i}$ and $W_{\bf i}^*$ as in \eqref{Whdef}, we obtain
\beas
\E[f'(W)]&=& \frac{1}{\sigma^2}\E\left(\sum_{\bf i} \sigma_{\bf i}^2 f'(W) + \sum_{{\bf i} \ne {\bf j}} \sigma_{{\bf i},{\bf j}} f'(W) \right)\\
         &=& \frac{1}{\sigma^2}\E\left(\sum_{\bf i} \xi_{\bf i}^2 f'(W) 
				 + \sum_{\bf i}\sum_{\substack{{\bf j}:{\bf j} \ne {\bf i}\\ |{\bf j}-{\bf i}|_\infty \le 1}}\xi_{\bf i}\xi_{\bf j}f'(W) 
	+ \sum_{\bf i}\sum_{\substack{{\bf j}:{\bf j} \ne {\bf i}\\ |{\bf j}-{\bf i}|_\infty > 1}}\sigma_{{\bf i},{\bf j}}f'(W) \right. \\
 &&\left.  + \sum_{\bf i} (\sigma_{\bf i}^2-\xi_{\bf i}^2) f'(W) 
+ \sum_{\bf i}\sum_{\substack{{\bf j}:{\bf j} \ne {\bf i}\\ |{\bf j}-{\bf i}|_\infty \le 1}} (\sigma_{{\bf i},{\bf j}}-\xi_{\bf i}\xi_{\bf j}) f'(W)\right) \\
         &=& \frac{1}{\sigma} \E \sum_{\bf i} \xi_{\bf i}W_{\bf i} f'(W)  
	+ \frac{1}{\sigma^2}\E\left( \sum_{\bf i}\sum_{\substack{{\bf j}:{\bf j} \ne {\bf i}\\ |{\bf j}-{\bf i}|_\infty > 1}}\sigma_{{\bf i},{\bf j}}f'(W) \right. \\
 &&\left.  + \sum_{\bf i} (\sigma_{\bf i}^2-\xi_{\bf i}^2) f'(W) 
+ \sum_{\bf i}\sum_{\substack{{\bf j}:{\bf j} \ne {\bf i}\\ |{\bf j}-{\bf i}|_\infty \le 1}} (\sigma_{{\bf i},{\bf j}}-\xi_{\bf i}\xi_{\bf j}) f'(W)\right)
\enas 
and
\beas
\E[Wf(W)] &=& \frac{1}{\sigma}\E\sum_{\bf i} \xi_{\bf i} f(W) =  \frac{1}{\sigma}\E\sum_{\bf i} \xi_{\bf i} f(W_{\bf i}^* + W_{\bf i}) \\
         &=& \frac{1}{\sigma}\E\sum_{\bf i} \left[ \xi_{\bf i} f(W_{\bf i}^*) + \xi_{\bf i}W_{\bf i} \int_0^1 f'(W_{\bf i}^* + u W_{\bf i}) du \right].
\enas
Substituting the two equations above into \eqref{steineq} with $w$ replaced by $W$, we have
\bea \label{proofeq}
\E[h(W)-N h]&=&\E[f'(W)-Wf(W)] \nn \\
&=& \E \left( \frac{1}{\sigma}\sum_{\bf i}\xi_{\bf i}W_{\bf i} \left(\int_0^1 \left(f'(W) -f'(W_{\bf i}^*+ uW_{\bf i})\right)du\right) \right. \nn \\
 &&+\frac{1}{\sigma^2}\sum_{\bf i}\sum_{\substack{{\bf j}:{\bf j} \ne {\bf i}\\ |{\bf j}-{\bf i}|_\infty > 1}}\sigma_{{\bf i},{\bf j}}f'(W) +  \frac{1}{\sigma^2}\sum_{\bf i} (\sigma_{\bf i}^2 - \xi_{\bf i}^2) f'(W) \nn \\
&& \left. + \frac{1}{\sigma^2}\sum_{\bf i}\sum_{\substack{{\bf j}:{\bf j} \ne {\bf i}\\ |{\bf j}-{\bf i}|_\infty \le 1}} (\sigma_{{\bf i},{\bf j}}-\xi_{\bf i}\xi_{\bf j}) f'(W) 
		   - \frac{1}{\sigma} \sum_{\bf i}  \xi_{\bf i} f(W_{\bf i}^*) \right).
\ena

Next we handle the five terms in \eqref{proofeq} separately. For the first term, we have
\bea \label{1st}
&&\frac{1}{\sigma}\left| \E  \sum_{\bf i}\xi_{\bf i} W_{\bf i} \int_0^1 \left(f'(W) -f'(W_{\bf i}^*+ uW_{\bf i})\right)du \right| \nn  \\
&& \hspace{15pt}= \frac{1}{\sigma}\left| \E \sum_{\bf i} \xi_{\bf i}W_{\bf i} \int_0^1 \int_{uW_{\bf i}}^{W_{\bf i}} f''(W_{\bf i}^*+t) dt du  \right|
\le \frac{2}{\sigma} \E \sum_{\bf i} |\xi_{\bf i}||W_{\bf i}| \left(  \int_0^1 \int_{u|W_{\bf i}|}^{|W_{\bf i}|} dt du\right) \nn \\
&& \hspace{15pt}=  \frac{1}{\sigma} \sum_{\bf i} \E |\xi_{\bf i}|W_{\bf i}^2 
   \le \frac{1}{\sigma} \sum_{\bf i} \left(\E |\xi_{\bf i}|^3\right)^{1/3}\left(\E |W_{\bf i}|^3\right)^{2/3}  
	= \frac{1}{\sigma}\sum_{\bf i} \left\|\xi_{\bf i}\right\|_3\left\|W_{\bf i}\right\|_3^2 \nn\\
&& \hspace{15pt}\le \frac{1}{\sigma^3} \sum_{\bf i}  \left(\sum_{{\bf s} \in D_{\bf i}^l} \left\|Y_{\bf s}\right\|_3\right) \left(\sum_{{\bf j}:|{\bf j}-{\bf i}|_\infty \le 1}\sum_{{\bf t} \in D_{\bf j}^l} \left\|Y_{\bf t}\right\|_3\right)^2 \nn  \\
&& \hspace{15pt} \le \frac{9^dm^dl^{3d}M^3}{\gamma^{3/2}n^{3d/2}} \le \frac{18^d M^3l^{2d}}{\gamma^{3/2}n^{d/2}},
\ena
where we have used the Holder's inequality in the second inequality, the triangle inequality in the third inequality, conditions \ref{b} and \ref{d} in the second last inequality and that $m \le 2n/l$ in the last inequality.

Moving on to the second term of \eqref{proofeq}, using the triangle inequality, we have
\beas
\sum_{\bf i}\sum_{\substack{{\bf j}:{\bf j} \ne {\bf i}\\ |{\bf j}-{\bf i}|_\infty > 1}}\left|\Cov(\xi_{\bf i},\xi_{\bf j})\right| \le 
  \sum_{{\bf i}\ne{\bf j}}\left|\Cov(\xi_{\bf i},\xi_{\bf j})\right| 
\le\sum_{{\bf i}\ne{\bf j}} \sum_{{\bf s} \in D_{\bf i}}\sum_{{\bf t} \in D_{\bf j}}\left|\Cov(Y_{\bf s},Y_{\bf t})\right| .
\enas
Applying Lemma \ref{lemrio} with $r=p=q=3$ and Lemma \ref{lempdl} using that $Y_{\bf i}$ is measurable with respect to $\mathcal{E}(C_{\bf i})$ and that $|C_{\bf i}| = 1$, the last expression is bounded by
\beas
&& \sum_{{\bf i}\ne{\bf j}} \sum_{{\bf s} \in D_{\bf i}}\sum_{{\bf t} \in D_{\bf j}}\alpha^{1/3}\left(\mathcal{E}(C_{\bf s}),\mathcal{E}(C_{\bf t})\right) \left\|Y_{\bf s}\right\|_3\left\|Y_{\bf t}\right\|_3 \nn \\
&&\hspace{15pt}\le M^2\sum_{{\bf i}\ne{\bf j}} \sum_{{\bf s} \in D_{\bf i}}\sum_{{\bf t} \in D_{\bf j}}\alpha^{1/3}_{1,1}\left(\dist(C_{\bf s},C_{\bf t})\right) \nn \\
&&\hspace{15pt}\le M^2\sum_{{\bf i}\ne{\bf j}} \sum_{{\bf s} \in D_{\bf i}}\sum_{{\bf t} \in D_{\bf j}}\sup_{|{\bf x}-{\bf y}|_\infty\ge \dist(C_{\bf s},C_{\bf t})}D({\bf x},{\bf y})^{1/3}.
\enas
Now using condition \ref{c}, we bound the last line by
\beas
&& \kappa^{1/3} M^2 \sum_{{\bf i}\ne{\bf j}} 
  \sum_{a_1,\ldots,a_d = -l+1}^{l-1}(l-|a_1|)\cdots(l-|a_d|)\exp\left(-\frac{\lambda}{3}\left|\begin{bmatrix} a_1+(j_1-i_1)l\\  \vdots \\ a_d+(j_d-i_d)l  \end{bmatrix}\right|_1 + \frac{d\lambda}{3}\right) \nn\\
	&& \hspace{15pt} = \kappa^{1/3} M^2 e^{d\lambda/3} \sum_{{\bf i}\ne{\bf j}}  \sum_{a_1,\ldots,a_d = -l+1}^{l-1}
	  \prod_{q=1}^d (l-|a_q|)\exp\left(-\frac{\lambda}{3} \left|a_q+(j_q-i_q)l\right|\right) \nn \\
&& \hspace{15pt} = \kappa^{1/3} M^2 e^{d\lambda/3}\sum_{{\bf i}\ne{\bf j}} \prod_{q=1}^d \sum_{a_q = -l+1}^{l-1}
	   (l-|a_q|)\exp\left(-\frac{\lambda}{3} \left|a_q+(j_q-i_q)l\right|\right)
\enas
where the factor of $e^{d\lambda/3}$ has arisen by the fact that each element-wise distance between $C_{\bf s}$ and $C_{\bf t}$ is $\max\{0,|a_k+(j_k-i_k)l|-1\}$ which is greater than or equal to $|a_k+(j_k-i_k)l|-1$.

Next, following the same argument as in equations (39)-(40) in the proof of Lemma 2.4 of \cite{GW18} with $\lambda$ replaced by $\lambda/3$, the last expression is equal to
\beas
&& \sum_{s=1}^d\kappa^{1/3}M^2 e^{d\lambda/3} {d \choose s} 2^s \sum_{\substack{1 \le i_k = j_k \le m \\ k=s+1,\ldots,d}}   
          \sum_{\substack{1 \le i_k < j_k \le m \\ k=1,\ldots,s}} \prod_{q=1}^s e^{-\frac{\lambda}{3}(j_q-i_q)l}\\
          && \text{\hspace{100pt}} \times\left(\prod_{q=1}^s \sum_{a_q=-l+1}^{l-1} (l-|a_q|) e^{-\frac{\lambda}{3} a_q}\right) 
					                                                   \left(\prod_{q=s+1}^d \sum_{a_q=-l+1}^{l-1} (l-|a_q|)e^{-\frac{\lambda}{3} |a_q|} \right)\\
&=& \sum_{s=1}^d\kappa^{1/3}M^2 e^{d\lambda/3}  {d \choose s} 2^s  m^{d-s}  \left(\sum_{k=1}^{m-1}(m-k)e^{-\frac{\lambda}{3} k l}\right)^s\\
          && \text{\hspace{25pt}} \times\left(l+\sum_{a=1}^{l-1} (l-a) \left(e^{\frac{\lambda}{3} a}+e^{-\frac{\lambda}{3} a}\right)\right)^s\left(l+2\sum_{b=1}^{l-1} (l-b)e^{-\frac{\lambda}{3} b}\right)^{d-s}.
\enas 
Invoking the identities in \eqref{Sum.krk4} and \eqref{Sum.krk3}, following the same computations as in (40) of \cite{GW18}, the last expression can be bounded by
\beas
\frac{\kappa^{1/3}M^2  n^d\left((4\mu_{\lambda}+2\nu_{\lambda})^d-\left(2\nu_{\lambda}\right)^d\right)}{l} ,
\enas 
where $\mu_\lambda$ and $\nu_\lambda$ are given as in \eqref{munudef}. Using \eqref{steineq:bounds} and condition \ref{d}, the second term is bounded by
\bea \label{2nd}
 \frac{1}{\sigma^2}\left|\sum_{\bf i}\sum_{\substack{{\bf j}:{\bf j} \ne {\bf i}\\ |{\bf j}-{\bf i}|_\infty > 1}}\sigma_{{\bf i},{\bf j}}f'(W)\right| 
    &\le& \frac{\sqrt{2}\kappa^{1/3}M^2 \left((4\mu_{\lambda}+2\nu_{\lambda})^d-\left(2\nu_{\lambda}\right)^d\right)}{\sqrt{\pi}\gamma l}.
\ena 
From this point, for simplicity, we will not mention when and how we use the assumptions \ref{a}, \ref{b}, \ref{c} and \ref{d} since they will be applied exactly as used in the first two terms. 

For the third term, we first bound it by 3 terms: $I_1$, $I_2$ and $I_3$ as follows,
\beas
 &&\frac{1}{\sigma^2}\left|\E \sum_{\bf i} (\sigma_{\bf i}^2 - \xi_{\bf i}^2) f'(W)\right| \nn\\
 &&  \hspace{15pt} = \frac{1}{\sigma^2} \left|\E \sum_{\bf i}(\sigma_{\bf i}^2 - \xi_{\bf i}^2)f'(W_{\bf i}^*) + \sum_{\bf i} (\sigma_{\bf i}^2 - \xi_{\bf i}^2)\left(f'(W)-f'(W_{\bf i}^*)\right)\right| \nn \\
&& \hspace{15pt} \le \frac{1}{\sigma^2} \sum_{\bf i}\left|\Cov\left( \xi_{\bf i}^2,f'(W_{\bf i}^*)\right)\right| + \frac{1}{\sigma^2} \sum_{\bf i} \sigma_{\bf i}^2 \E\left|f'(W)-f'(W_{\bf i}^*)\right| \nn\\
&& \hspace{50pt}  +   \frac{1}{\sigma^2}\sum_{\bf i} \E \xi_{\bf i}^2\left|f'(W)-f'(W_{\bf i}^*)\right| := I_1+I_2+I_3. \nn\\
\enas
For $I_1$, we have
\bea \label{3rdi1}
I_1 &\le& \frac{1}{\sigma^2}\sum_{\bf i} \alpha^{1/3}\left(\mathcal{E}\left(\bigcup_{{\bf t} \in D_{\bf i}}C_{\bf t}\right),\mathcal{E}\left(\bigcup_{{\bf j}:|{\bf j}-{\bf i}|_\infty > 1}\bigcup_{{\bf t} \in D_{\bf j}}C_{\bf t}\right)\right) \left\|\xi_{\bf i}\right\|_3^2 \left\|f'(W_{\bf i}^*)\right\|_\infty \nn\\
    &\le& \frac{1}{\sigma^2} \sqrt{\frac{2}{\pi}}l^{2d}M^2 \sum_{\bf i}\alpha^{1/3}_{l^d,n^d-3^d l^d}(l) \nn\\
		&\le& \frac{1}{\sigma^2} \sqrt{\frac{2}{\pi}}l^{7d/3}m^d(n^d-3^dl^d)^{1/3}M^2 \sup_{|{\bf x}-{\bf y}|\ge l} |D({\bf x},{\bf y})|^{1/3} \nn \\
		&\le& \frac{2^{d+1/2}\kappa^{1/3} M^2 l^{4d/3} n^{d/3}}{\sqrt{\pi}\gamma e^{\lambda l/3}},
\ena
where we have used Lemma \ref{lemrio} with $r=3,p=3/2,q=\infty$ and that for any random variable $X$ that the third moment exists, $\left\|X^2\right\|_{3/2} = \left\|X\right\|_3^2$ and $\left\|X\right\|_p \le \left\|X\right\|_q$ for $p \le q$ in the first inequality, the absolute bound in \eqref{steineq:bounds} and the facts that $\dist\left(\bigcup_{{\bf t} \in D_{\bf i}}C_{\bf t},\bigcup_{{\bf j}:|{\bf j}-{\bf i}|_\infty > 1}\bigcup_{{\bf t} \in D_{\bf j}}C_{\bf t}\right) \ge l$ and $|C_{\bf i}|=1$ in the second inequality and Lemma \ref{lempdl} in the third inequality.
For $I_2$, again using the absolute bound in \eqref{steineq:bounds} and that  $\left\|X\right\|_p \le \left\|X\right\|_q$ for $p \le q$ we obtain
\bea \label{3rdi2}
I_2 &\le&  \frac{1}{\sigma^2} \sum_{\bf i} \sigma_{\bf i}^2 \left|f''\right|_{\infty}\E|W-W_{\bf i}^*| 
    \le \frac{2}{\sigma^2} \sum_{\bf i} \E \xi_{\bf i}^2 \E|W_{\bf i}|\nn \\ 
		&\le& \frac{2}{\sigma^2} \sum_{\bf i} \left\|\xi_{\bf i}\right\|_3^2 \left\|W_{\bf i}\right\|_3 \le \frac{2\cdot 3^d m^dl^{3d} M^3}{\gamma^{3/2} n^{3d/2}}
		      \le \frac{2^{d+1} 3^d  M^3 l^{2d}}{\gamma^{3/2} n^{d/2}}.
\ena
For $I_3$, using the Holder's inequality, we have
\bea \label{3rdi3}
I_3 &\le& \frac{2}{\sigma^2} \sum_{\bf i} \E \xi_{\bf i}^2|W_{\bf i}| \le \frac{2}{\sigma^2}\sum_{\bf i} \left\|\xi_{\bf i}\right\|_3^2\left\|W_{\bf i}\right\|_3 \nn \\
    &\le& \frac{2^{d+1} 3^d  M^3 l^{2d}}{\gamma^{3/2} n^{d/2}}.
\ena

For the fourth term, we again bound it by 3 terms: $J_1$, $J_2$ and $J_3$ as follows,
\beas
&& \frac{1}{\sigma^2}\left|\E\sum_{\bf i}\sum_{\substack{{\bf j}:{\bf j} \ne {\bf i}\\ |{\bf j}-{\bf i}|_\infty \le 1}} (\sigma_{{\bf i},{\bf j}}-\xi_{\bf i}\xi_{\bf j}) f'(W)\right| \nn \\
&& \hspace{15pt} =  \frac{1}{\sigma^2}\left|\E\sum_{\bf i}\sum_{\substack{{\bf j}:{\bf j} \ne {\bf i}\\ |{\bf j}-{\bf i}|_\infty \le 1}} \left[(\sigma_{{\bf i},{\bf j}}-\xi_{\bf i}\xi_{\bf j}) f'(W_{{\bf i},{\bf j}}^*) + (\sigma_{{\bf i},{\bf j}}-\xi_{\bf i}\xi_{\bf j})\left(f'(W)-f'(W_{{\bf i},{\bf j}}^*)\right)\right]\right| \nn\\
&& \hspace{15pt} \le \frac{1}{\sigma^2} \sum_{\bf i}\sum_{\substack{{\bf j}:{\bf j} \ne {\bf i}\\ |{\bf j}-{\bf i}|_\infty \le 1}} \left|\Cov\left(\xi_{\bf i}\xi_{\bf j},f'(W_{{\bf i},{\bf j}}^*)\right)\right| 
     + \frac{1}{\sigma^2} \sum_{\bf i}\sum_{\substack{{\bf j}:{\bf j} \ne {\bf i}\\ |{\bf j}-{\bf i}|_\infty \le 1}} |\sigma_{{\bf i},{\bf j}}|\E\left|f'(W)-f'(W_{{\bf i},{\bf j}}^*)\right|  \nn \\
	&&	 \hspace{80pt} + \frac{1}{\sigma^2} \sum_{\bf i}\sum_{\substack{{\bf j}:{\bf j} \ne {\bf i}\\ |{\bf j}-{\bf i}|_\infty \le 1}} \E |\xi_{\bf i}\xi_{\bf j}|  \left|f'(W)-f'(W_{{\bf i},{\bf j}}^*)\right| := J_1+J_2+J_3.
\enas
Proceeding similarly to the third term, we have 
\bea \label{4thj1}
J_1 &\le& \frac{1}{\sigma^2}\sum_{\bf i}\sum_{\substack{{\bf j}:{\bf j} \ne {\bf i}\\ |{\bf j}-{\bf i}|_\infty \le 1}} \alpha^{1/3}\left(\mathcal{E}\left(\bigcup_{{\bf t} \in D_{\bf i}\cup D_{\bf j}}C_{\bf t}\right),\mathcal{E}\left(\bigcup_{\substack{{\bf k}:|{\bf k}-{\bf i}|_\infty > 1\\ |{\bf k}-{\bf j}|_\infty > 1}}\bigcup_{{\bf t} \in D_{\bf k}}C_{\bf t}\right)\right) \left\|\xi_{\bf i}\right\|_3 \left\|\xi_{\bf j}\right\|_3  \left\|f'(W_{{\bf i},{\bf j}}^*)\right\|_\infty \nn\\
    &\le& \frac{1}{\sigma^2} \sqrt{\frac{2}{\pi}}l^{2d}M^2 \sum_{\bf i}\sum_{\substack{{\bf j}:{\bf j} \ne {\bf i}\\ |{\bf j}-{\bf i}|_\infty \le 1}}\alpha^{1/3}_{2l^d,n^d-4\cdot 3^{d-1}l^d}(l) \nn\\
		&\le& \frac{2^{1/3}}{\sigma^2} \sqrt{\frac{2}{\pi}} (3^d-1)l^{7d/3}m^d(n^d-4\cdot 3^{d-1}l^d)^{1/3}M^2 \sup_{|{\bf x}-{\bf y}|\ge l} |D({\bf x},{\bf y})|^{1/3} \nn \\
		&\le& \frac{2^{d+5/6}3^d \kappa^{1/3} l^{4d/3}n^{4d/3}M^2}{\sqrt{\pi}\gamma n^d e^{\lambda l/3}}
		\le \frac{2\cdot 6^d \kappa^{1/3} M^2l^{4d/3}n^{d/3}}{\sqrt{\pi} \gamma  e^{\lambda l/3}},
\ena
\bea \label{4thj2}
J_2 &\le&  \frac{1}{\sigma^2} \sum_{\bf i}\sum_{\substack{{\bf j}:{\bf j} \ne {\bf i}\\ |{\bf j}-{\bf i}|_\infty \le 1}}  |\sigma_{{\bf i},{\bf j}}| \left|f''\right|_{\infty}\E|W-W_{{\bf i},{\bf j}}^*| 
    \le \frac{2}{\sigma^2} \sum_{\bf i}\sum_{\substack{{\bf j}:{\bf j} \ne {\bf i}\\ |{\bf j}-{\bf i}|_\infty \le 1}}  \E |\xi_{\bf i}\xi_{\bf j}| \E|W_{{\bf i},{\bf j}}|\nn \\ 
		&\le& \frac{2}{\sigma^2} \sum_{\bf i}\sum_{\substack{{\bf j}:{\bf j} \ne {\bf i}\\ |{\bf j}-{\bf i}|_\infty \le 1}} \left\|\xi_{\bf i}\right\|_3 \left\|\xi_{\bf j}\right\|_3\left\|W_{{\bf i},{\bf j}}\right\|_3 \le \frac{8\cdot 3^{2d-1} m^dl^{3d} M^3}{\gamma^{3/2} n^{3d/2}}
		      \le \frac{2^{d+3} 3^{2d-1} M^3 l^{2d}}{\gamma^{3/2} n^{d/2}}
\ena
and
\bea \label{4thj3}
J_3 &\le& \frac{2}{\sigma^2} \sum_{\bf i}\sum_{\substack{{\bf j}:{\bf j} \ne {\bf i}\\ |{\bf j}-{\bf i}|_\infty \le 1}}\E |\xi_{\bf i}\xi_{\bf j} W_{{\bf i},{\bf j}}| \le \frac{2}{\sigma^2}\sum_{\bf i} \left\|\xi_{\bf i}\xi_{\bf j}\right\|_{3/2}\left\|W_{{\bf i},{\bf j}}\right\|_3 \nn \\
    &\le& \frac{2}{\sigma^2} \sum_{\bf i}\sum_{\substack{{\bf j}:{\bf j} \ne {\bf i}\\ |{\bf j}-{\bf i}|_\infty \le 1}} \left\|\xi_{\bf i}\right\|_3 \left\|\xi_{\bf j}\right\|_3\left\|W_{{\bf i},{\bf j}}\right\|_3
		    \le \frac{2^{d+3} 3^{2d-1} M^3 l^{2d}}{\gamma^{3/2} n^{d/2}}.
\ena

Now we handle the last term. Using Lemma \ref{lemrio} with $r=3/2,p=3,q=\infty$ and Lemma \ref{lempdl} along with the facts used in the first four terms, we obtain

\bea \label{5th}
\frac{1}{\sigma} \left|\E \sum_{\bf i}  \xi_{\bf i} f(W_{\bf i}^*)\right| &=& \frac{1}{\sigma} \sum_{\bf i}  \left|\Cov\left(\xi_{\bf i},f(W_{\bf i}^*)\right)\right| \nn\\
&\le& \frac{1}{\sigma} \sum_{\bf i} \sum_{{\bf t} \in D_{\bf i}} \left|\Cov\left(Y_{\bf t},f(W_{\bf i}^*)\right)\right| \nn \\
&\le& \frac{1}{\sigma} \sum_{\bf i} \sum_{{\bf t} \in D_{\bf i}} \alpha^{2/3}\left(\mathcal{E}(C_{\bf t}),\mathcal{E}\left(\bigcup_{{\bf j}:|{\bf j}-{\bf i}|_\infty > 1}\bigcup_{{\bf s} \in D_{\bf j}}C_{\bf s}\right)\right) \left\|Y_{\bf t}\right\|_3 \left\|f(W_{\bf i}^*)\right\|_\infty \nn \\
&\le& \frac{2M}{\sigma} \sum_{\bf i} \sum_{{\bf t} \in D_{\bf i}} \alpha_{1,n^d-3^dl^d}^{2/3}(l) \nn \\
&\le& \frac{2M m^d l^d (n^d-3^dl^d)^{2/3}}{\sigma}\sup_{|{\bf x}-{\bf y}|\ge l} |D({\bf x},{\bf y})|^{2/3} \nn \\
&\le& \frac{2M \kappa^{2/3} m^d l^d (n^d-3^dl^d)^{2/3}}{\gamma^{1/2} n^{d/2} e^{2\lambda l/3}} \le \frac{2^{d+1} M \kappa^{2/3} n^{7d/6}}{\gamma^{1/2}  e^{2 \lambda l/3}}.
\ena

Combining \eqref{1st} to \eqref{5th}, we have
\bea \label{combine}
&& \frac{18^d M^3l^{2d}}{\gamma^{3/2}n^{d/2}}+\frac{\sqrt{2}\kappa^{1/3}M^2 \left((4\mu_{\lambda}+2\nu_{\lambda})^d-\left(2\nu_{\lambda}\right)^d\right)}{\sqrt{\pi}\gamma l}+\frac{2^{d+1/2}\kappa^{1/3} M^2 l^{4d/3} n^{d/3}}{\sqrt{\pi}\gamma e^{\lambda l/3}} + \frac{2^{d+2} 3^d  M^3 l^{2d}}{\gamma^{3/2} n^{d/2}} \nn \\
&& \hspace{20pt}+ \frac{2\cdot 6^d \kappa^{1/3} M^2l^{4d/3}n^{d/3}}{\sqrt{\pi} \gamma  e^{\lambda l/3}} +\frac{2^{d+4} 3^{2d-1} M^3 l^{2d}}{\gamma^{3/2} n^{d/2}} + \frac{2^{d+1} M \kappa^{2/3} n^{7d/6}}{\gamma^{1/2}  e^{2\lambda l/3}} \nn \\
&& = \frac{9\cdot 18^d M^3l^{2d}}{\gamma^{3/2}n^{d/2}}+\frac{\sqrt{2}\kappa^{1/3}M^2 \left((4\mu_{\lambda}+2\nu_{\lambda})^d-\left(2\nu_{\lambda}\right)^d\right)}{\sqrt{\pi}\gamma l} \nn \\
&& \hspace{20pt} + \frac{3\cdot 6^d\kappa^{1/3} M^2l^{4d/3}n^{d/3}}{\sqrt{\pi} \gamma  e^{\lambda l/3}}
   +\frac{2^{d+1} M \kappa^{2/3} n^{7d/6}}{\gamma^{1/2}  e^{2\lambda l/3}} . 
\ena

From the last expression, it is clear that if $l$ is of order $n^s$ for some $s>0$ then the last two terms do not contribute to the rate as $l$ is the exponents in the denominators of both terms. Therefore, we just seek for $s>0$ that $l = n^s$ provides the best rate of convergence for the first two terms in \eqref{combine}. By using the technique in (44) of \cite{GW18} that the minimum of $al^{2d}+b/l$ is achieved at $l_0 = (b/2ad)^{1/{(2d+1)}}$ and taking $l=\lfloor{l_0}\rfloor$, it follows that 
\beas
l = \left\lfloor{\left( \frac{\sqrt{2\gamma}\kappa^{1/3} \left((4\mu_{\lambda}+2\nu_{\lambda})^d-\left(2\nu_{\lambda}\right)^d\right)n^{d/2}}{18^{d+1}  \sqrt{\pi} M }\right)^{1/(2d+1)}}\right\rfloor
\enas
and hence
\beas
d_1\left(\mathcal{L}(S_{{\bf k}}^n/\sigma_{n,{\bf k}}),\mathcal{L}(Z)\right) \le 
\frac{C_{1,d,M,\kappa,\gamma}}{n^{d/(4d+2)}} + \frac{C_{2,d,M,\kappa,\gamma} n^{d(4d+1)/(6d+3)}}{\exp\left(\theta_{d,M,\kappa,\gamma}n^{d/(4d+2)}\right)}
       + \frac{C_{3,d,M,\kappa,\gamma} n^{7d/6}}{\exp\left(2\theta_{d,M,\kappa,\gamma}n^{d/(4d+2)}\right)},
\enas
where $\theta_{d,M,\kappa,\gamma}$ is defined in \eqref{thetadef} and $C_{1,d,M,\kappa,\gamma}$, $C_{2,d,M,\kappa,\gamma}$ and $C_{3,d,M,\kappa,\gamma}$ are given in \eqref{Cdef}.

\bbox

\section{Applications to determinantal point processes} \label{sec:app}

In this section, we obtain $L^1$ bounds between functionals of determinantal point processes on $\mathbb{R}^d$ and the normal. Before moving on to our goal, we first provide the definition and an existence condition for determinantal point processes.

\begin{definition} \label{DPPs}
Let $\mathbb{K}: \left(\mathbb{R}^d\right)^2 \rightarrow \mathbb{C}$ be a measurable function. A point process $X$ on $\mathbb{R}^d$ is said to be a determinantal point process with kernel $\mathbb{K}$ if it is simple and its joint intensities with respect to the Lebesgue measure satisfy
\beas
\rho_n({\bf x}_1,\ldots,{\bf x}_n) = \det\left[\mathbb{K}({\bf x}_i,{\bf x}_j)\right]_{1 \le i,j \le n},
\enas
for every $n \in \mathbb{N}_1$ and ${\bf x}_1,\ldots,{\bf x}_n \in \mathbb{R}^d$.
\end{definition} 

By the definition and \eqref{ddef}, it follows immediately that $D({\bf x},{\bf y}) = -\mathbb{K}^2({\bf x},{\bf y})$. Due to this fact, in Theorem \ref{mainthm2} below we assume \eqref{as1} on $\mathbb{K}^2({\bf x},{\bf y})$ instead of $D({\bf x},{\bf y})$.
Recall that, for a kernel $\mathbb{K}$ that is locally square integrable, an associated integral operator is defined as 
\beas
\mathcal{K}f({\bf x}) = \int_{\mathbb{R}^d} \mathbb{K}({\bf x},{\bf y})f({\bf y}) d {\bf y}  \text{ \ \ for a.e. } {\bf x} \in \mathbb{R}^d,
\enas
for functions $f \in L^2(\mathbb{R}^d)$ that vanish a.e. outside a bounded subset of $\mathbb{R}^d$. For a compact set $S \subset \mathbb{R}^d$, the restriction of $\mathcal{K}$ to $S$, denoted by $\mathcal{K}_S$, is the bounded linear operator on $L^2(S)$ defined by
\beas
\mathcal{K}_Sf({\bf x}) = \int_{S} \mathbb{K}({\bf x},{\bf y})f({\bf y}) d {\bf y}  \text{ \ \ for a.e. } {\bf x} \in S.
\enas 
We say that $\mathcal{K}_S$ is of \textit{trace class} if $\sum_j |\lambda_j^S| < \infty$ where $\lambda_j^S$ are eigenvalues of $\mathcal{K}_S$ and $\mathcal{K}$ is \textit{locally of trace class} if $\mathcal{K}_S$ is of trace class for all compact subsets $S \subset \mathbb{R}^d$.

Next we present an existence condition from \cite{Mac75}, proved probabilistically in \cite{HKPV09}. Note that the results are actually in a more general space.

\begin{theorem} [\cite{Mac75}, \cite{HKPV09}]
Let $\mathbb{K}: \left(\mathbb{R}^d\right)^2 \rightarrow \mathbb{C}$ be locally square integrable Hermitian measurable function and its associated integral operator $\mathcal{K}$ is locally of trace class. Then a determinantal point process $X$ with kernel $\mathbb{K}$ exists if and only if all eigenvalues of $\mathcal{K}$ are contained in $[0,1]$.
\end{theorem}

Therefore, in the following, we consider only $\mathbb{K}$ that is locally square integrable and Hermitian such that its associated operator is locally of trace class with eigenvalues in $[0,1]$.

In this section, for $X$ a determinantal point process on $\mathbb{R}^d$ with $d \in \mathbb{N}_1$, we consider a function $f: \Omega \rightarrow \mathbb{R}$ defined by
\bea \label{fdef}
f(Y) = \sum_{S \subset Y} g(S),
\ena
where $g$ is a bounded function such that $g(S)=0$ when $\diam(S)> \tau$ for some fixed $\tau>0$. Here we denote $\diam(S)=\sup_{x,y \in S} \left|x-y\right|_\infty$. Our interest is in obtaining $L^1$ bounds to the normal for the distribution of $f(X \cap \Lambda_n)$ where $\Lambda_n = [0,n]^d$.

A CLT for the number of points in $X \cap \Lambda_n$ or $N(\Lambda_n)$ which corresponds to \eqref{fdef} with $g(S) = \mathbf{1}_{|S|=1}$ was proved by \cite{Sos00a} under some weak restriction on its variance decay that even allows a logarithmic rate. The result was generalized to spacing variables in \cite{Sos00b}. A CLT for the general functionals of the form \eqref{fdef} was shown in \cite{PDL17} with a stronger restriction on the variance decaying rate along with a few more assumptions.

The following theorem obtains an $L^1$ bound between the standardized $f(X \cap \Lambda_n)$ and the normal. Note that we will not track the constant since it requires complicated calculations and we use a couple of inequalities from \cite{PDL17} that do not provide explicit constants. Hence $C$ in the proof that follows may change from line to line.

\begin{theorem} \label{mainthm2}
Let $n \in \mathbb{N}_1$ and $X$ be a determinantal point process with kernel $\mathbb{K}$ that is bounded, locally square integrable and Hermitian such that its associated operator is locally of trace class with eigenvalues in $[0,1]$. Let $f$ be defined as in \eqref{fdef} with $g$ be bounded and $g(S)=0$ when $\diam(S)> \tau$ for some fixed $\tau>0$. Letting $\sigma_n^2 = \Var(f(X \cap \Lambda_n))$ with $\Lambda_n = [0,n]^d$, assume that 
\bea \label{as1}
\sup_{|{\bf x}-{\bf y}|_\infty \ge r} \mathbb{K}({\bf x},{\bf y}) \le O(e^{-\lambda r}) \text{ \ \ for some \ } \lambda >0,  \text{ \ \ and}
\ena
\bea \label{as2}
\liminf_n \sigma_n^2/n^d >0.
\ena
Then, with $W_n = (f(X \cap \Lambda_n)-\E f(X \cap \Lambda_n))/\sigma_n$, there exists $C>0$ such that
\bea \label{mainbound2}
d_1\left(\mathcal{L}(W_n),\mathcal{L}(Z)\right) \le \frac{C}{n^{d/(4d+2)}}, 
\ena
where $Z$ is a standard normal random variable.
\end{theorem}


\proof
We follow the same technique as in the proof of Theorem 4.4 of \cite{PDL17} approximating $f(X \cap \Lambda_n)$ by 
\beas
f_{\tilde{\Lambda}_n}(X) = \sum_{S \subset X} g(S)\mathbf{1}_{\tilde{\Lambda}_n}(S^0),
\enas
where $\tilde{\Lambda}_n = \bigcup_{{\bf i} \in I_n}C_{{\bf i}}$, $I_n = \{{\bf i}: C_{\bf i}\oplus\tau \subset \Lambda_n \}$ and $S^0$ is the barycenter of the set $S$. By the definition, it is clear that
\beas
f_{\tilde{\Lambda}_n}(X)  = \sum_{{\bf i} \in I_n} f_{C_{\bf i}}(X) .
\enas

It follows from \cite{PDL17} that 
\bea \label{vardiffbound}
\Var\left(f(X \cap \Lambda_n)-f_{\tilde{\Lambda}_n}(X)\right) &=& O(|\Lambda_n\backslash \tilde{\Lambda}_n\oplus \tau|) \nn\\
   &\le& C n^{d-1} ,
\ena
and
\bea \label{covdiffbound}
\left|\Cov\left(f(X \cap \Lambda_n),f(X \cap \Lambda_n)-f_{\tilde{\Lambda}_n}(X)\right)\right| &\le& \sigma_n\sqrt{\Var\left(f(X \cap \Lambda_n)-f_{\tilde{\Lambda}_n}(X)\right)} \nn\\ 
	&\le& C\sigma_n n^{(d-1)/2}.
\ena
Denoting $\tilde{\sigma}_n = \Var(f_{\tilde{\Lambda}_n}(X))$, we have
\bea \label{vararg}
\tilde{\sigma}_n^2 &=& \Var(f_{\tilde{\Lambda}_n}(X)-f(X \cap \Lambda_n)+f(X \cap \Lambda_n)) \nn\\
   &=& \sigma_n^2 + \Var\left(f(X \cap \Lambda_n)-f_{\tilde{\Lambda}_n}(X)\right) \nn\\
	   && \hspace{10pt}- \Cov\left(f(X \cap \Lambda_n),f(X \cap \Lambda_n)-f_{\tilde{\Lambda}_n}(X)\right) .
\ena
Using \eqref{vardiffbound} and \eqref{covdiffbound} and the assumption \eqref{as2}, we have 
\beas
|\tilde{\sigma}_n^2-\sigma_n^2| \le Cn^{d-1/2},
\enas
and hence,
\bea \label{diffvar}
|\tilde{\sigma}_n-\sigma_n| \le \frac{Cn^{d-1/2}}{\tilde{\sigma}_n+\sigma_n} \le \frac{Cn^{d-1/2}}{\sigma_n} \le Cn^{(d-1)/2}.
\ena

Using the triangle inequality, we have
\beas
d_1\left(\mathcal{L}(W_n),\mathcal{L}(Z)\right) &\le& d_1\left(\frac{1}{\sigma_n}(f(X \cap \Lambda_n)-\E f(X \cap \Lambda_n)),\frac{1}{\sigma_n}(f_{\tilde{\Lambda}_n}(X) -\E f_{\tilde{\Lambda}_n}(X) )\right) \\
   && + d_1\left(\frac{1}{\sigma_n}(f_{\tilde{\Lambda}_n}(X) -\E f_{\tilde{\Lambda}_n}(X) ),\frac{1}{\tilde{\sigma}_n}(f_{\tilde{\Lambda}_n}(X) -\E f_{\tilde{\Lambda}_n}(X) )\right) \\
	&& +d_1\left(\frac{1}{\tilde{\sigma}_n}(f_{\tilde{\Lambda}_n}(X) -\E f_{\tilde{\Lambda}_n}(X) ),Z\right) := D_1+D_2+D_3.
\enas
We handle $D_1$, $D_2$ and $D_3$ separately. Using the definition of the $L^1$ distance in \eqref{Wasdef} and applying \eqref{vardiffbound}, \eqref{covdiffbound} and \eqref{as2} for the first two terms, we obtain
\beas
D_1 &\le& \frac{1}{\sigma_n}\E|f(X \cap \Lambda_n)-f_{\tilde{\Lambda}_n}(X) -\E(f(X \cap \Lambda_n)-f_{\tilde{\Lambda}_n}(X))| \\
    &\le& \frac{1}{\sigma_n}\sqrt{\Var(f(X \cap \Lambda_n)-f_{\tilde{\Lambda}_n}(X))} \le \frac{Cn^{(d-1)/2}}{\sigma_n} \le Cn^{-1/2},
\enas
and
\beas
D_2 &\le& \left|\frac{1}{\sigma_n}-\frac{1}{\tilde{\sigma}_n}\right|\E|f_{\tilde{\Lambda}_n}(X) -\E f_{\tilde{\Lambda}_n}(X)| \\
    &\le& \left|\frac{1}{\sigma_n}-\frac{1}{\tilde{\sigma}_n}\right|\tilde{\sigma}_n = \frac{|\tilde{\sigma}_n-\sigma_n|}{\sigma_n}\le \frac{Cn^{(d-1)/2}}{\sigma_n} \le Cn^{-1/2},
\enas
where we have also used \eqref{diffvar} in the second last inequality of $D_2$.

Finally, we bound the last term by applying Theorem \ref{mainthm1} showing that conditions \ref{a}-\ref{d} there are satisfied. Condition \ref{a} is obvious by the definition of determinantal point processes, \ref{b} follows immediately from \cite{PDL17}, \ref{c} follows from the fact that $D({\bf x},{\bf y}) = -\mathbb{K}^2({\bf x},{\bf y})$ and the assumption \eqref{as1}, and \ref{d} follows from \eqref{vararg} and the assumption \eqref{as2}. Therefore, $D_3 \le C n^{-d/(4d+2)}$ which proves the claim since $D_1$ and $D_2$ are at least of order $n^{-1/2}$ and thus converge to zero faster than $n^{-d/(4d+2)}$.

\bbox

Next we state a few remarks and a corollary that relate to Theorem \ref{mainthm2}.

\begin{remark} \label{varsuffcond}
If $g$ in \eqref{fdef} has support on the set $\{S \subset \mathbb{R}^d:|S|=p\}$ for some $p \in \mathbb{N}_1$ then a sufficient condition for $\Var(f(X \cap \Lambda_n))$ to satisfy the assumption \eqref{as2} was given in Lemma B.7 of \cite{PDL17}. That is, \eqref{as2} holds if $\left\|\mathcal{K}\right\|<1$ and
\bea \label{condb7}
\liminf_n \frac{1}{n^d}\int_{\Lambda_n^p}g(\{{\bf x}_1,\ldots,{\bf x}_p\}) \det[\mathbb{K}({\bf x}_i,{\bf x}_j)]_{1\le i,j\le p} d({\bf x}_1, \ldots, {\bf x}_p)>0.
\ena
\end{remark}

Recall that the Laguerre-Gaussian family of determinantal point processes (\cite{BL16}) is a family of determinantal point processes. As the name may suggest, point processes in this family have Laguerre-Gaussian functions as their kernel functions, that is, for ${\bf x},{\bf y} \in \mathbb{R}^d$, $\mathbb{K}({\bf x},{\bf y}) = \mathbb{C}_{m,\alpha,\rho}({\bf x}-{\bf y})$ where for $m \in \mathbb{N}_1$, $\alpha> 0$, $\rho > 0$ and ${\bf z} \in \mathbb{R}^d$,
\bea \label{cdef}
\mathbb{C}_{m,\alpha,\rho}({\bf z}) = \frac{\rho}{{m-1+d/2 \choose m-1}}L_{m-1}^{d/2}\left(\frac{1}{m}\left|\frac{{\bf z}}{\alpha}\right|^2\right)e^{-|{\bf z}/\alpha|^2/m},
\ena 
and
\bea \label{Ldef}
L_n^s(x) = \sum_{k=0}^n {n+s \choose n-k} \frac{(-x)^k}{k!} \text{ \ for all $x,s \in \mathbb{R}$ \ and $n \in \mathbb{N}_0$}
\ena
be the Laguerre polynomials. When $m=1$, it is known as Gaussian determenantal point process which is one of the most well known processes in the Laguerre-Gaussian family. It was used earlier to analyze data sets in spatial statistics (see \cite{LMR15} for example).

It is clear by the definition that point processes in this family satisfy the assumption \eqref{as1}. The following corollary is a special case of Theorem \ref{mainthm2} when $X$ is in this family. The result follows from Lemma B.7 of \cite{PDL17} and the fact shown in \cite{BL16} that $\left\|\mathcal{K}\right\|<1$ when \eqref{alphacond} holds with strict inequalities.

\begin{corollary}
Let $X$ be a determenantal point process in the Laguerre-Gaussian family with 
\bea \label{alphacond}
0 < \alpha \le \left[{m-1+d/2 \choose m-1}/\rho(m\pi)^{d/2}\right]^{1/d}.
\ena
Let $f$ and $g$ be defined as in the statement of Theorem \ref{mainthm2}. Assume that \eqref{as2} is satisfied, then the bound \eqref{mainbound2} in Theorem \ref{mainthm2} holds. Furthermore, if $g$ in \eqref{fdef} has support on the set $\{S \subset \mathbb{R}^d:|S|=p\}$ for some $p \in \mathbb{N}_1$, $\alpha$ satisfies \eqref{alphacond} with strict inequalities and \eqref{condb7} is satisfied then the bound \eqref{mainbound2} holds.
\end{corollary}

A point process is said to be \textit{translation-invariant} or \textit{stationary} if its kernel $\mathbb{K}({\bf x},{\bf y}) = \mathbb{C}({\bf x}-{\bf y})$ for some function $\mathbb{C}:\mathbb{R}^d \rightarrow \mathbb{R}$. It is clear by the definition that determinantal point processes in the Laguerre-Gaussian family are stationary. The work \cite{Sos00b} calculated the variance of $N([-n,n]^d)$ for $X$ a stationary point process on $\mathbb{R}^d$ and provided an example that satisfies a weaker version of our assumption \eqref{as2}.

\begin{remark}
Let $X$ be a stationary point process on $\mathbb{R}$ with kernel $\mathbb{C}$ and the Fourier transform of $\mathbb{C}$, $
\hat{\mathbb{C}}(x) = \mathbf{1}_B(x)$ with $B=\bigsqcup_{i\ge 1}[i,i+1/i^\gamma]$ for some $\gamma>1$. By \cite{Sos00b}, $\Var(N([-n,n])) = O(n^{1/\gamma})$ for $n \in \mathbb{N}_1$. $N([n,n])$ can be represented by $f(X \cap [-n,n])$ where f is as in \eqref{fdef} with $g(S) = \mathbf{1}_{|S|=1}$. As mentioned in Remark \ref{nsnd}, our method also works when $\Var(f(X \cap [-n,n])) = O(n^s)$ with $s>6/7$ which corresponds to this example when $1< \gamma < 7/6$. However, the assumption \eqref{as1} is not met here.
\end{remark}

\begin{remark}
The Gaussian unitary ensemble is one of the most well-known determinantal point processes on $\mathbb{R}$ (see \cite{HKPV09}) with Kernel
\beas
\mathbb{K}_n(x,y) = \sum_{k=0}^{n-1}H_k(x)H_k(y)
\enas
with respect to the Gaussian measure $d\mu(x) = \frac{1}{\sqrt{2\pi}}e^{-x^2/2}dx$ where $H_k(\cdot)$ are Hermite polynomials, obtained by applying Gram-Schmidt orthogonalization procedure to $\{1,x,x^2,\ldots\}$ in $L^2(\mathbb{R},\mu)$. It is obvious that in this case the assumption \eqref{as1} is satisfied. However, the sufficient condition mentioned in Remark \ref{varsuffcond} does not apply to this process since $\left\|\mathcal{K}\right\|=1$. 
\end{remark}



\section*{Acknowledgements}

The author would like to thank an anonymous reviewer for extremely helpful comments and detailed suggestions that lead to an improvement of the organization of the paper. The author would like to also thank the Institute for Mathematical Sciences, National University of Singapore for the hospitality during the time the author edited the second version of this paper.

\end{document}